\documentclass[10pt]{asme2ej}
\usepackage{cite}
\usepackage{times}
\usepackage{makeidx}
\usepackage{amsmath}
\usepackage{graphics}
\usepackage{graphicx}
\usepackage{subfigure}
\usepackage{amssymb}
\usepackage{euler}
\usepackage{supertabular}

\title{Probabilistic Control for Uncertain Systems}

\author{Randa Herzallah
\affiliation{FET \\ Al-Balqa' Applied University, Jordan \\Email:herzallah.r@gmail.com
}}

\begin{document}

\maketitle


\begin{abstract}
{\it
In this paper a new framework has been applied to the design of
controllers which encompasses nonlinearity, hysteresis and
arbitrary density functions of forward models and inverse
controllers. Using mixture density networks, the probabilistic
models of both the forward and inverse dynamics are estimated such
that they are dependent on the state and the control input. The
optimal control strategy is then derived which minimizes
uncertainty of the closed loop system. In the absence of reliable
plant models, the proposed control algorithm incorporates
uncertainties in model parameters, observations, and latent
processes. The local stability of the closed loop system has been
established. The efficacy of the control algorithm is demonstrated
on two nonlinear stochastic control examples with additive and
multiplicative noise.}

\end{abstract}

\section{Introduction}\label{sec:Intro}
The conventional theory of stochastic control is particularly
suitable for taking into account randomly varying system
parameters and designing probabilistic control strategies under
uncertain working conditions. However, the functional equations
describing its solution are mostly computationally infeasible and
are solved approximately by assuming the certainty equivalence
principle, for example, and the ignorance of model
uncertainty~\cite{Astrom89,Narendra94,Yaz86}.

Machine learning has been proposed in~\cite{Herzallah08} to model
the conditional distribution of the system dynamics. The
application of these methods to control systems is still in its
infancy where only the Gaussian distribution of the process noise
dynamics has been estimated~\cite{Herzallah08,Herzallah07c}.
Nevertheless, these methods are proven to be promising in
addressing a number of key weaknesses in traditional control
approaches that are either based on deterministic models, ignore
model uncertainty or treat uncertainty as a nuisance
parameter. A key objective of this paper is to extend these
currently well developed methods of estimating the conditional
distributions of the systems dynamics in such a way that they are
able to approximate the general distributions of the dynamical
stochastic and deterministic nonlinear systems and exploit them in
control strategies. This will be accomplished by using the mixture density network (MDN) from the neural
network field~\cite{Bishop95}. In its original structure the MDN has
been proposed to model the probability density functions of stochastic static models. In previous work~\cite{Herzallah04}, we
have successfully extended the idea of the MDN such that the
probability density functions of the stochastic dynamic models can
be estimated as well. The developed theory was then used to
estimate the conditional distribution of an ill posed inverse
controller where a unique solution cannot be found. In this work the
MDN will be used to estimate the general distributions of system
models which are not constrained by the Gaussian assumption.
Although, the idea of a mixture density network is not
new~\cite{Evans00,Richmond03,Bishop95,Herzallah03a}, no results to
our knowledge on using this neural network  have been reported in
the control literature.

Several attempts to deriving control strategies of stochastic
uncertain systems have been discussed in the
literature~\cite{Hayakawa09,Pin09,Primbs09,Mirkin08,Petersen08,Zhang09}.
However, most of the current
methods~\cite{Primbs09,Mirkin08,Petersen08,Zhang09} are devoted to linear uncertain systems, and do not consider, to a large
extent, functional uncertainties as a result of unknown latent
dynamics and hysteresis characteristics of the process. The closed loop entropy has been proposed in~\cite{Yue03} to characterize the uncertainty of the tracking error for general nonlinear and non--Gaussian stochastic systems.  However the probability density function (pdf) of the tracking
error in~\cite{Yue03} is measured by assuming a known pdf  of the
random input that affects the dynamic of the system. In other
words, the system output is assumed to be invertible with respect
to the noise. This is a crude assumption for many general
practical systems since it is often difficult to measure the pdf
of arbitrary random inputs. It also implicitly indicates the
existence of an accurate model that describes the system dynamics,
and consequently ignores model uncertainty.

In this paper a new framework which does not assume the
invertibility of the system output with respect to the noise or
even the invertibility of the control input with respect to the
system output is developed. It also does not assume the existence
of accurate models that describe the system forward and inverse
dynamics. The new framework is based on estimating the
probabilistic models of the tracking error and the inverse
controller from process data. These probabilistic models are not
constrained by Gaussian assumptions. They are state and control
input dependent and are estimated using a MDN. The control input is then designed to minimize uncertainty of the closed loop system. Our framework provides a theoretical but practically implementable control mechanism that leads to a more robust and efficient control strategy under highly complex working conditions. The innovation of this paper stems from accepting uncertainty as fundamental to the understanding of the control problem, placing it at the heart of a probabilistic generative view of control.
This is counter to the traditional view of control in which,
when uncertainty is considered at all, it is simply modeled and
taken to be a separable nuisance aspect of the more basic
deterministic control strategy. Moreover, the method developed in
this paper is suitable for controlling multi--modal control
systems, where a unique solution does not exist.

\section{Statement of the problem}\label{sec:Stat}
The general stochastic control problem is considered for systems
with input $u_k$, a measurable state vector ${\bf z}_k =
[\mathrm{y}_{k-1}, \\ \mathrm{y}_{k-2}, \ldots, \mathrm{y}_{k-n},
u_{k-1}, \ldots, u_{k-m}]$, and future values of the system
output, $\mathrm{y}_{k+\tau}$, affected by a random force
$\varepsilon_{k+\tau}$. The system behavior can be represented by
the following nonlinear ARMAX stochastic model
\begin{equation}
\mathrm{y}_{k+\tau} = f({\bf z}_k , u_k, \varepsilon_{k+\tau})
\label{eq:StochasticEq1},
\end{equation}
where $\tau$ is the relative degree of the system, and $f(.)$ is an
unknown nonlinear function that represents the system dynamic. In general $f$ need not be invertible with respect to the random input, $\varepsilon_{k+\tau}$. Besides,
there is no assumption made on whether $\varepsilon_{k+\tau}$ has
a known pdf or is an independent and identically distributed
random process. Because of the existence of the random force, only
the probability distribution of the future output values, $P_{\mathrm{y}}$  can be specified from the state and control at each instant of time, $k$. The aim of control is subject to some constraints such as a finite energy budget,
\begin{eqnarray}
e_{k+\tau} &=& f({\bf z}_k , u_k, \varepsilon_{k+\tau}) -
\mathrm{y}^d_{k+\tau}, \nonumber \\
&=& g(\mathrm{y}_{k+\tau}, \mathrm{y}^d_{k+\tau}),
\label{eq:StochasticEq3}
\end{eqnarray}
where $g(.)$ is the stochastic model of the tracking error and is
obtained by subtracting the desired output $\mathrm{y}^d_{k+\tau}$
from the system stochastic function $f({\bf z}_k , u_k,
\varepsilon_{k+\tau})$. Hence $g(.)$ also need not be invertible with respect to $\varepsilon_{k+\tau}$. Moreover, the tracking error distribution is not necessarily Gaussian because $\varepsilon_{k+\tau}$ is a general stochastic noise. The density of $e_{k+\tau}$ can be obtained from the density of $\mathrm{y}_{k+\tau}$ as
\begin{equation}
P_e(\mathrm{y}_{k+\tau}, \mathrm{y}^d_{k+\tau}) =
P_{\mathrm{y}}\bigg(e_{k+\tau} + \mathrm{y}^d_{k+\tau} \bigg).
\label{eq:StochasticEq5}
\end{equation}
In this formulation the value of the variable $u_k$ parameterizes
and affects the distribution of the system output and consequently
the distribution of the tracking error. For the system to perform
well in practice, the controller should be designed such that the
pdf of the tracking error is made as narrow as possible. A narrow
distribution indicates that the uncertainty of the tracking error
is small which also corresponds to a small variance. In addition,
the mean value of the tracking error and control energy should
also be minimized. For this purpose, the performance index is
set to be,
\begin{eqnarray}
J(u_k) = R \int \parallel e_{k+\tau}-<e_{k+\tau}\mid
\mathrm{y}_{k+\tau}, \mathrm{y}^d_{k+\tau}> \parallel^2
P_e(\mathrm{y}_{k+\tau}, \mathrm{y}^d_{k+\tau}) d e_{k+\tau}
+ M \big( \int e_{k+\tau} P_e(\mathrm{y}_{k+\tau},
\mathrm{y}^d_{k+\tau}) d e_{k+\tau} \big)^2 + Q u^2_k,
\label{eq:StochasticEq6}
\end{eqnarray}
where $R, M$ and $Q$ are constant weights for the variance, mean,
and control input respectively.

\textit{\textbf{Remark 1:}} To reemphasize, the pdf of the tracking
error need not be invertible with respect to the random force $\epsilon_{k+\tau}$, but is invertible with respect to the system output as
specified by~\eqref{eq:StochasticEq5}. Therefore, both the
output pdf and the tracking error pdf will be used in this article
mutually, and they have the same parameters with the mean of the
tracking error pdf being different by the desired output than the
mean of the output pdf. This will be further discussed in
Section~\ref{sec:TrackErr}

Since the pdf of the system output is not constrained by the
Gaussian assumption, a sum of squares or cross entropy error
function for estimating the system output are not expected to
yield satisfactory results. For this purpose a new class of
network models obtained by combining a conventional neural network
with a mixture of Gaussians is proposed in this paper to estimate
the conditional probability distribution of the system output from
process data. It is called a mixture density network and can in
principle represent arbitrary conditional probability
distributions in the same way that a conventional neural network can
represent arbitrary functions. The control objective is then to
track the specified desired output, $\mathrm{y}^d_{k+\tau}$ to a
small neighborhood of zero with the output
$\mathrm{y}_{k+\tau}$, while ensuring local stability of the
system output. The problem therefore arises as to how inverse
controllers are acquired from the general distribution of the
system output. Such learning must be able to divide up the control
into appropriate regions which can be recombined to generate
the system behavior. In this paper we propose a novel framework which
can solve the learning of the general distributions of the system
output and inverse controller in a computationally coherent manner
from a single principle. The basic idea of the proposed framework
is that a mixture of Gaussians exist to control the system and each
is augmented with a corresponding Gaussian of the system output.
It consists of a model pool of couples of MDNs of inverse
controllers and system models. Each couple evaluates a number of
probable control signals, and the couple generating the most
suitable control signal is used to control the system. This framework is especially useful for estimating the general distributions of systems from process data. As well be demonstrated shortly, it is efficient for controlling complex systems with large uncertainties which are not constrained by the Gaussian assumption. Besides, it gives superior results for controlling systems characterized by hysteresis and multi--modality.

\subsection{Mixture density networks and tracking error distribution}\label{sec:TrackErr}
\begin{figure*}[htbp]
\centerline{\subfigure[]{\includegraphics[width=2.5in]{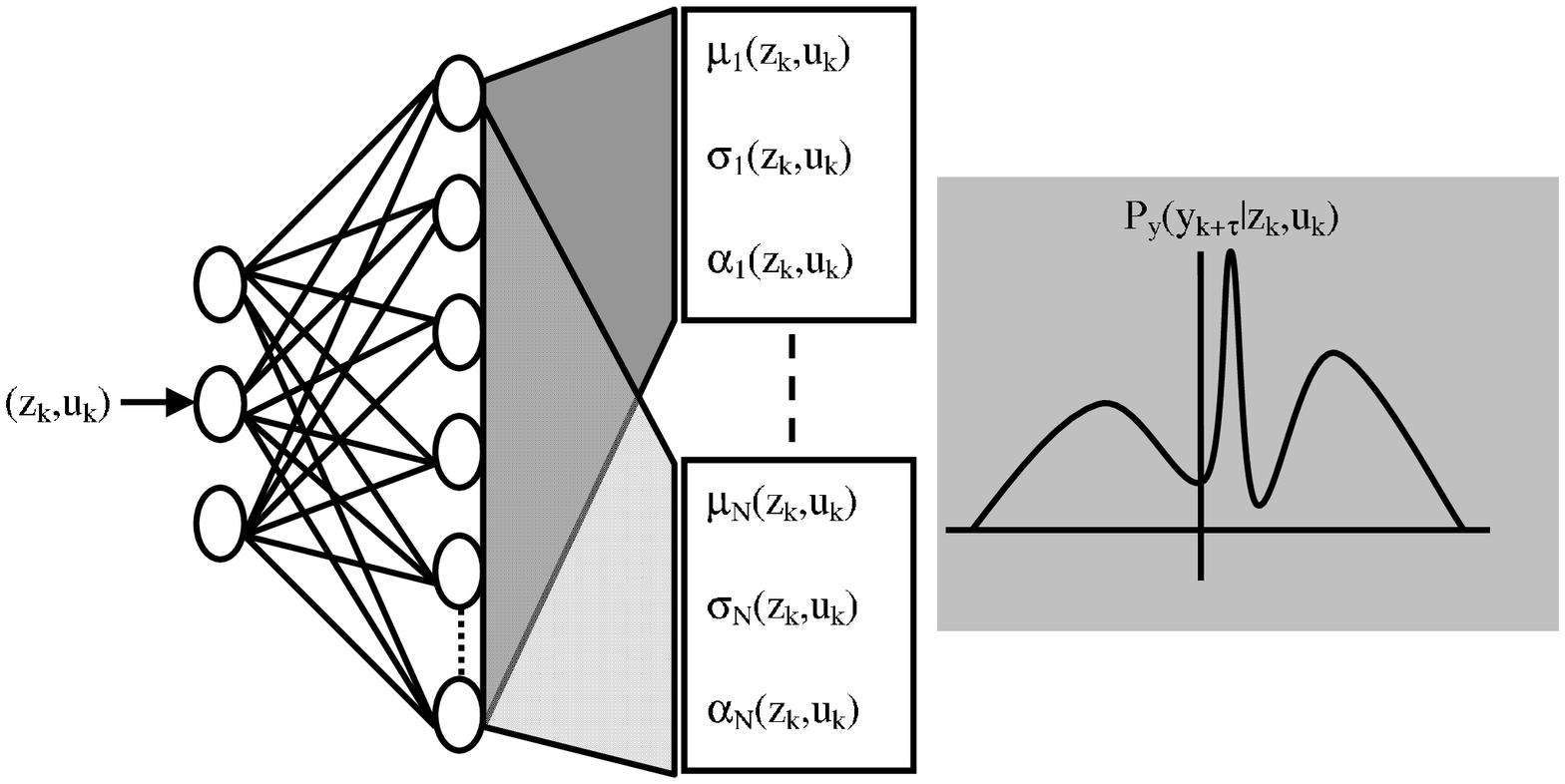}
\label{fig_first_case}} \hfil
\subfigure[]{\includegraphics[width=2.5in]{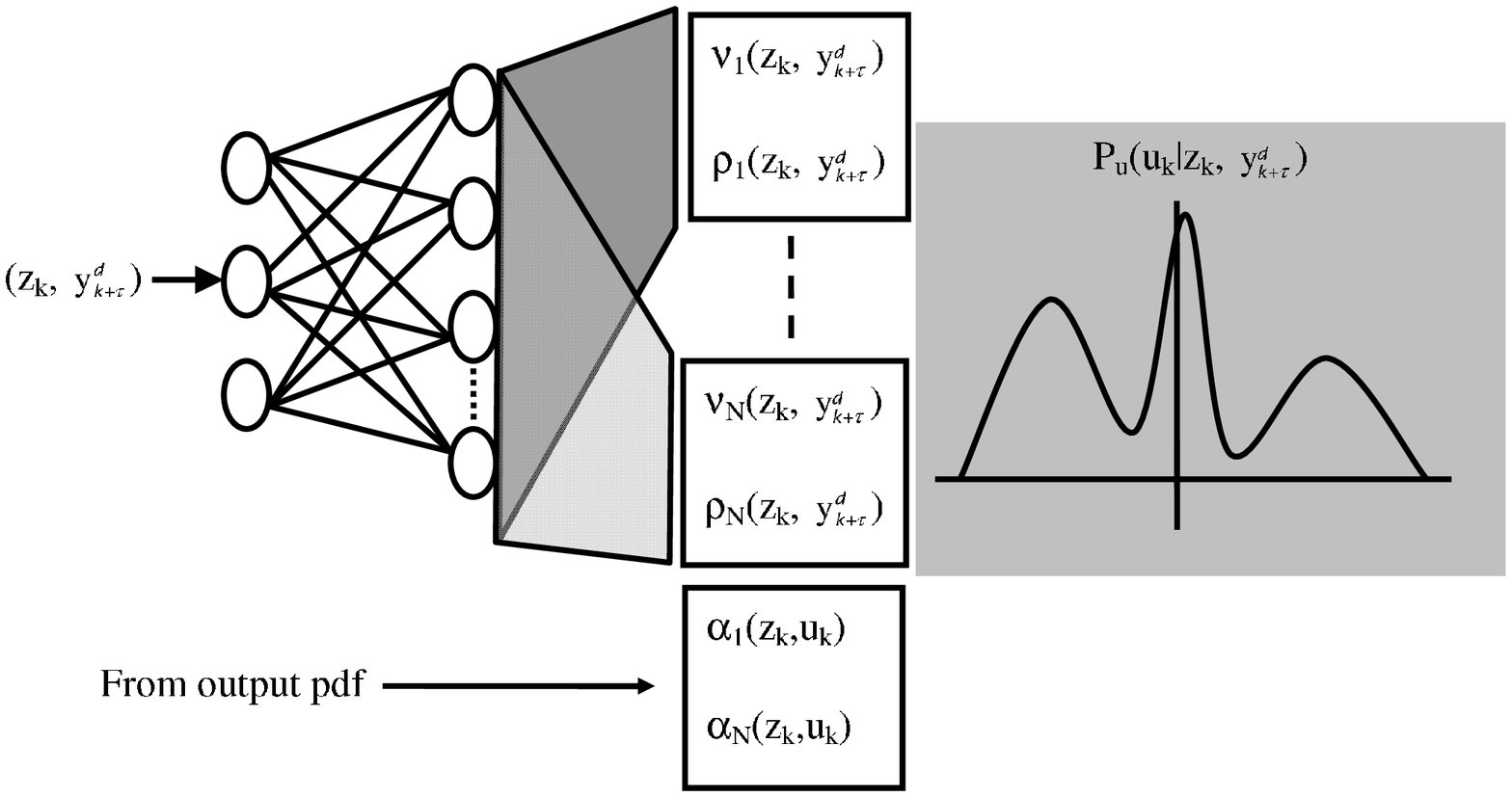}
\label{fig_second_case}}} \caption{Architecture of forward and inverse MDNs: (a) The architecture of the system output MDN. (b) The architecture of the inverse controller MDN.} \label{fig:ArchMix}
\end{figure*}

Mixture density networks have been employed in
many system identification and inverse
applications~\cite{Kravchenko09,Evans00} and have also been shown to
provide a general framework for approximating the conditional
distribution of the inverse controller where multi--modality and
hysteresis play critical roles~\cite{Herzallah04}. In a mixture
density network, the probability density of the system output is
given by a linear combination of Kernel functions in the form
\begin{equation}
P_{\mathrm{y}}(\mathrm{y}_{k+\tau} \mid {\bf z}_k , u_k, \theta) =
\sum_{i=1}^N \alpha_i({\bf z}_k , u_k) \phi_i (\mathrm{y}_{k+\tau}
\mid {\bf z}_k , u_k), \label{eq:StochasticEq7}
\end{equation}
where $N$ is the number of kernels in the mixture. The
parameters $\alpha_i({\bf z}_k, u_k)$ are the prior probabilities
of $\mathrm{y}_{k+\tau}$ having been generated from the $i^{th}$
component of the mixture. The functions $\phi_i
(\mathrm{y}_{k+\tau} \mid {\bf z}_k , u_k)$ represent the
conditional density of $\mathrm{y}_{k+\tau}$ for the $i^{th}$
kernel. Various choices of the kernel functions are
possible~\cite{Bishop95}, however, in this article Gaussian kernel
functions are considered,
\begin{equation}
\phi_i (\mathrm{y}_{k+\tau} \mid {\bf z}_k , u_k) =
\frac{1}{\sqrt{2 \pi} \sigma^2_i({\bf z}_k , u_k)} \exp \bigg \{ -
\frac{\parallel \mathrm{y}_{k+\tau} - \mu_i ({\bf z}_k , u_k)
\parallel^2}{2 \sigma^2_i({\bf z}_k , u_k)} \bigg \}, \label{eq:StochasticEq8}
\end{equation}
where $\mu_i ({\bf z}_k , u_k)$, and $\sigma^2_i({\bf z}_k , u_k)$
represent the centre and the variance respectively of the $i^{th}$
kernel. Note that the prior probabilities, the centre, and the
variances of the output pdf are taken to be continuous functions
of the input variables $({\bf z}_k , u_k)$. These functions are
estimated as the outputs of a feed--forward neural network that
takes $({\bf z}_k , u_k)$ as input. They represent the set of
parameters which govern the system output distribution, and are
denoted by $\theta$ in Equation~\eqref{eq:StochasticEq7}. This
combination of a density model and a feed--forward neural network
is represented schematically in Figure~\ref{fig:ArchMix}, a. Gaussian kernels as specified by Equation~\eqref{eq:StochasticEq8} can approximate any given density function to arbitrary accuracy.

Using $\mathrm{y}_{k+\tau}$ as a target, the mixture density
network is then trained to minimize the negative logarithm of the
probability density function of the system output by using
back-propagation
\begin{equation}
E = -\ln \bigg \{ \sum_{i=1}^N \alpha_i({\bf z}_k , u_k) \phi_i
(\mathrm{y}_{k+\tau} \mid {\bf z}_k , u_k) \bigg \}.
\label{eq:StochasticEq8p}
\end{equation}
Details of the derivatives of the error function
\eqref{eq:StochasticEq8p} with respect to the outputs of the
networks and constraints of the mixture models can be found
in~\cite{Herzallah04,Bishop95}.

Using Equations~\eqref{eq:StochasticEq5},
~\eqref{eq:StochasticEq7} and~\eqref{eq:StochasticEq8}, the
distribution of the tracking error can then be obtained from the
estimated conditional distribution of the system output as follows
\begin{equation}
P_e(e_{k+\tau} \mid \mathrm{y}_{k+\tau},\mathrm{y}^d_{k+\tau},
\theta) = \sum_{i=1}^N \alpha_i({\bf z}_k , u_k) \frac{1}{\sqrt{2
\pi} \sigma^2_i({\bf z}_k , u_k)} \exp \bigg \{ - \frac{\parallel
\mathrm{e}_{k+\tau} - \ell_i ({\bf z}_k , u_k,
\mathrm{y}^d_{k+\tau})
\parallel^2}{2 \sigma^2_i({\bf z}_k , u_k)} \bigg \}, \label{eq:StochasticEqet}
\end{equation}
where $\ell_i ({\bf z}_k , u_k, \mathrm{y}^d_{k+\tau}) = \mu_i
({\bf z}_k , u_k) - \mathrm{y}^d_{k+\tau}$ is the centre of the
pdf of the tracking error which is obtained by subtracting the
desired output value from the centre of the output pdf.

\subsection{Mixture density networks and the distribution of the inverse controller}\label{sec:Dist}
Based on the ability of the multiple Gaussians to describe the
density function of the system output, we suggest that for each
behavior captured by a kernel function, it is desired to learn a
control strategy or in other words a paired inverse kernel
function should be designed. As such, the probability density
function of the inverse controller can be generated as the
summation of outputs from these inverse kernels weighted by the
prior probabilities of the output pdf,
\begin{equation}
P_u (u_k \mid {\bf z}_k , \mathrm{y}^d_{k+\tau}, \vartheta) =
\sum_{i=1}^N \alpha_i({\bf z}_k , u_k) \psi_i (u_k \mid {\bf z}_k
, \mathrm{y}^d_{k+\tau}), \label{eq:StochasticEq9}
\end{equation}
where $\psi_i (u_k \mid {\bf z}_k , \mathrm{y}^d_{k+\tau})$
represent the conditional density of $u_k$ for the $i^{th}$
kernel.

\textbf{\textit{Remark 2:}} Note that in
Equation~\eqref{eq:StochasticEq9} the prior probabilities of the
inverse controller are fixed to those obtained from the output
pdf. The aim is that each inverse kernel learns to provide
a suitable control signal under the context for which its paired
kernel of the output pdf most likely produces the output value.

The kernel functions of the inverse controller are again taken to
be Gaussian kernel functions,
\begin{equation}
\psi_i (u_k \mid {\bf z}_k , \mathrm{y}^d_{k+\tau}) =
\frac{1}{\sqrt{2 \pi} \rho^2_i({\bf z}_k , \mathrm{y}^d_{k+\tau})}
\exp \bigg \{ - \frac{\parallel u_k - \nu_i ({\bf z}_k ,
\mathrm{y}^d_{k+\tau})
\parallel^2}{2 \rho^2_i({\bf z}_k , \mathrm{y}^d_{k+\tau})} \bigg
\}, \label{eq:StochasticEq10}
\end{equation}
where $\nu_i ({\bf z}_k , \mathrm{y}^d_{k+\tau})$, and
$\rho^2_i({\bf z}_k , \mathrm{y}^d_{k+\tau})$ represent the centre
and the variance respectively of the $i^{th}$ kernel of the
inverse controller. Similar to the output pdf, the centre, and the
variances of the inverse controller pdf are continuous functions
of the input variables $({\bf z}_k , \mathrm{y}^d_{k+\tau})$. The
centre and the variances functions are estimated as the outputs of
a feed--forward neural network that takes $({\bf z}_k ,
\mathrm{y}^d_{k+\tau})$ as input. This combination of a density
model and a feed--forward neural network is represented
schematically in Figure~\ref{fig:ArchMix}, b. The centre, variances
and priors of the inverse controller MDN represent the set of
parameters which govern the inverse controller distribution, and
are denoted by $\vartheta$ in Equation~\eqref{eq:StochasticEq9}. Here the target of the mixture density network is the optimal
control input as calculated in Section~\ref{sec:OpCont},
Equation~\eqref{eq:StochasticEq13} or its linearized form
specified in Equation~\eqref{eq:StochasticEq19}. The mixture
density network is then trained to minimize the negative logarithm
of the probability density function of the control input $u_k$ by using
back-propagation.

To emphasize, the output pdf takes the control signal and the
state values as inputs and estimates the conditional density
function of the output of the system. The controller pdf takes the
desired output of the system and the state values as inputs and
estimates the conditional density function of the control input.
The prior probabilities of the inverse controller are fixed to
their corresponding probabilities from the output pdf.
Conceptually speaking if the output of the system is most likely
produced by one of the kernels, its corresponding inverse kernel
receives the major part of the error signal and its output
contribute significantly to the conditional density of the control
input. Fixing the prior probabilities of the inverse controller to
those obtained from the output pdf, is realistic in practice. It
actually ensures that the output pdf and its counterpart of the
inverse controller are tightly coupled both through training and
control phases. Of course, one can still choose to have independent
priors for the inverse controller.

A key advantage of using the MDN is its ability to represent arbitrary
conditional probability distributions in the same way a
conventional neural network can represent arbitrary functions.
Moreover, the parameters of MDNs are optimized such as to minimize
the negative log-likelihood of the probability density
functions~\eqref{eq:StochasticEq7} and~\eqref{eq:StochasticEq9},
therefore, a complete description of the estimated output can be
obtained. This is of particular interest to control problems in
which the mapping to be learned is
multi--valued~\cite{Herzallah03a} as often arises in the solution
of inverse control problems and robotics
applications~\cite{White92,Molina04} and hysteretic nonlinear
systems~\cite{Ren09,Ikhouane08}. The above framework provides a new generic modeling and control methods appropriate for active control of complex uncertain systems seeking stochastically optimal control strategies in
systems which exhibit nonlinearity, hysteresis, multimodality,
randomness and uncertainty.
\subsection{Control Algorithm Design}\label{sec:OpCont}
The development so far assumes no prior information about known
pdfs of the system output or the inverse controller. All pdfs
required are estimated using MDNs from process data. The
parameters of those conditional distributions (means, variances
and priors) are continuous functions of the input variables. This
allows the development of a pragmatic method for estimating and
incorporating functional uncertainties in deriving the optimal
control law.

The probabilistic control problem is a nonlinear optimization
problem that can be solved by setting the derivative of the
performance function~\eqref{eq:StochasticEq6} with respect to the
control signal equal to zero,
\begin{equation}
\frac{\partial J(u_k)}{\partial u_k} = 0.
\label{eq:StochasticEq11}
\end{equation}
Using Equations~\eqref{eq:StochasticEqet},~\eqref{eq:StochasticEq6} and~\eqref{eq:StochasticEq11} gives,
\begin{eqnarray}
&R \bigg \{ \sum_i \frac{\partial \alpha_i({\bf z}_k ,
u_k)}{\partial u_k} \bigg \{ \sigma^2_i({\bf z}_k , u_k) + \bigg [
\ell_i ({\bf z}_k , u_k, \mathrm{y}^d_{k+\tau}) - \sum_l
\alpha_l({\bf z}_k , u_k) \ell_l ({\bf z}_k , u_k,
\mathrm{y}^d_{k+\tau}) \bigg ]^2 \bigg \} \nonumber \\ &+ \sum_i
\alpha_i({\bf z}_k , u_k) \bigg \{ \frac{\partial \sigma^2_i({\bf
z}_k , u_k)}{\partial u_k} + 2 \bigg [ \ell_i ({\bf z}_k , u_k,
\mathrm{y}^d_{k+\tau}) - \sum_l \alpha_l({\bf z}_k , u_k) \ell_l
({\bf z}_k , u_k, \mathrm{y}^d_{k+\tau}) \bigg ] \nonumber
\\ &\times \bigg [ \frac{\partial \ell_i ({\bf z}_k , u_k,
\mathrm{y}^d_{k+\tau})}{\partial u_k} - \sum_l \bigg (
\frac{\partial \alpha_l({\bf z}_k , u_k)}{\partial u_k} \ell_l
({\bf z}_k , u_k, \mathrm{y}^d_{k+\tau}) + \alpha_l({\bf z}_k ,
u_k) \frac{\partial \ell_l ({\bf z}_k , u_k,
\mathrm{y}^d_{k+\tau})}{\partial u_k} \bigg) \bigg] \bigg \}
\bigg\} \nonumber \\ &+ 2 M \sum_i \alpha_i({\bf z}_k , u_k) \ell_i ({\bf
z}_k , u_k, \mathrm{y}^d_{k+\tau}) \sum_i \bigg \{ \frac{\partial
\alpha_i({\bf z}_k , u_k)}{\partial u_k} \ell_i ({\bf z}_k , u_k,
\mathrm{y}^d_{k+\tau})+ \alpha_i({\bf z}_k , u_k) \frac{\partial
\ell_i ({\bf z}_k , u_k, \mathrm{y}^d_{k+\tau})}{\partial u_k}
\bigg \} + 2 Q u_k = 0 \nonumber \\ &\varphi(\chi_k,u_k)+ 2 Q u_k = 0, \label{eq:StochasticEq13}
\end{eqnarray}
where $\chi_k = (z_k,\mathrm{y}^d_{k+\tau}) \in R^{n+m+1}$. The solution to~\eqref{eq:StochasticEq13} cannot be analytically
obtained due to the nonlinear nature of the parameters of the
tracking error pdf. As such, this solution can only guarantee the
search for local minima. An alternative approach for
solving~\eqref{eq:StochasticEq13} analytically, would be to
formulate a recursive algorithm for the control input $u_k$.

Note that $\varphi(\chi_k,u_k)$ is continuous because the parameters of
the tracking error pdf are continuous and first order
differentiable with respect to $u_k$. Applying the first backward difference operator, $\bigtriangledown$ to~\eqref{eq:StochasticEq13} yields the
following recursive formula for $u_k$
\begin{equation}
\varphi(\chi_k,u_k) - \varphi(\chi_{k-1},u_{k-1}) = - 2 Q (u_k -
u_{k-1}). \label{eq:StochasticEq17}
\end{equation}
The increment of $\varphi(\chi_k,u_k)$ can be approximated from
the first order Taylor expansion of $\varphi(\chi_k,u_k)$ around
the previous operating point $(\chi_{k-1},u_{k-1})$:
\begin{equation}
\varphi(\chi_k,u_k) - \varphi(\chi_{k-1},u_{k-1}) \approx
\frac{\partial \varphi}{\partial \chi} \bigg |_{k-1} (\chi_k -
\chi_{k-1}) + \frac{\partial \varphi}{\partial u} \bigg |_{k-1}
(u_k - u_{k-1}). \label{eq:StochasticEq18}
\end{equation}
Substitute~\eqref{eq:StochasticEq18}
into~\eqref{eq:StochasticEq17}, it follows that the optimal
control input at time $k$ can be calculated using,
\begin{equation}
u^*_k = u^*_{k-1} - \frac{\frac{\partial \varphi}{\partial \chi}
\bigg |_{k-1}}{2Q+ \frac{\partial \varphi}{\partial u} \bigg
|_{k-1}} (\chi_k - \chi_{k-1}). \label{eq:StochasticEq19}
\end{equation}
The control input obtained in~\eqref{eq:StochasticEq19} can then
be used as the target value for estimating the conditional
distribution of the inverse controller as specified
by~\eqref{eq:StochasticEq9}.

Several methods for calculating the output of the mixture density
network have been proposed in the literature. In multimodal
control problems for example, the distribution of the system
output will consist of limited numbers of distinct branches. In
this case one specific branch from the estimated conditional
density of the MDN more likely is needed to be selected. Two
examples of how to select a specific branch are the most likely,
and the most probable output values. In this paper the most probable
output value corresponding to the most probable branch will be used. Since each component of the mixture model is normalized, the most probable branch is given by
\begin{eqnarray}
\mathrm{arg} \hspace{0.1cm} \underset{i}{\mathrm{max}} \{\alpha_i
\} \longrightarrow
u_k = \mu_i. \label{eq:StochasticEq20}
\end{eqnarray}
To summarize, the following algorithm can be readily obtained:
\begin{enumerate}
\item At sample time $k$, solve equation~\eqref{eq:StochasticEq13}
and calculate the optimal control input, or alternatively use
equation~\eqref{eq:StochasticEq19} to calculate the optimal
control input.
\item Update the parameters of the MDN that estimates the
conditional distribution of the inverse controller.
\item Calculate the control input $u_k$ from the controller MDN.
\item Apply $u_k$ to the stochastic
system~\eqref{eq:StochasticEq1}.
\item Update the parameters of the MDN that estimates the output
probability density function.
\item Obtain the density function of the tracking error
from~\eqref{eq:StochasticEqet}.
\end{enumerate}
This control algorithm can specifically be implemented on-line while actually performing in a good manner since knowledge of uncertainty is taken into consideration. The forward model of the plant to be controlled, and the controller can both be adapted on--line which means that speed requirement becomes stringent. However, since on--line optimisation is kept local fast convergence of the networks is expected. Moreover, the above control algorithm  provides a general solution for stochastic systems subject to arbitrary random inputs with unknown pdfs. Indeed it is a general solution for stochastic and deterministic systems characterized by functional uncertainty. No assumptions
are made about the invertibility of the system output with respect
to the random inputs or the invertibility of the controller with
respect to the system output. Thus, it can be concluded that this
control formulation can be widely applied to many general
practical stochastic non--Gaussian systems characterized by
hysteresis and high levels of complexity and uncertainty. The scientific quality of our proposed method over other recently
developed stochastic control methods focusses around the
acceptance of `noise' as `intrinsic uncertainty' which often
cannot be ignored and is absolutely key to solving the control
problem.
\section{Local Stability Analysis}
In this section, we consider the local stability of the closed
loop system given by~\eqref{eq:StochasticEq1}
and~\eqref{eq:StochasticEq20}.

\textbf{\textit{Theorem 1:}} Assume that the plant described
by~\eqref{eq:StochasticEq1} is being controlled
by~\eqref{eq:StochasticEq20}. Assume the following:
\begin{enumerate}
\item There exist a set of parameters $\theta$ of the tracking
error distribution such that $\mathcal{D}(P_e(\mathrm{y}_{k+\tau},
\mathrm{y}^d_{k+\tau}) \parallel P^t_e(\mathrm{y}_{k+\tau},
\mathrm{y}^d_{k+\tau})) \leq \epsilon_1$. Here $\mathcal{D}(.,.)$
denotes the Kullback Leibler divergence distance, $P^t_e(.)$
denotes the true or ideal distribution of the tracking error, and
$\epsilon_1$ is a positive constant which represents the
estimation error of the tracking error pdf.
\item There exist a set of parameters $\vartheta$ of the controller
distribution such that $\mathcal{D}(P_u({\bf z}_k,
\mathrm{y}^d_{k+\tau}) \parallel P^t_u({\bf z}_k,
\mathrm{y}^d_{k+\tau})) \leq \epsilon_2$. Here $P^t_u(.)$ denotes
the true or ideal distribution of the control input, and
$\epsilon_2$ is a positive number which represents the estimation
error of the pdf of the controller.
\end{enumerate}
Then, the closed loop system formed by~\eqref{eq:StochasticEq1}
and~\eqref{eq:StochasticEq20} is locally stable.

\textit{Proof:} To prove the local stability of the closed loop
system,~\eqref{eq:StochasticEq1} is linearized to give
\begin{equation}
\bigtriangleup \mathrm{y}_{k} = \sum_{i=1}^n \frac{\partial
f}{\partial \mathrm{y}_{k-i}} \bigtriangleup \mathrm{y}_{k-i} +
\sum_{j=0}^m \frac{\partial f}{\partial u_{k-j}} \bigtriangleup
u_{k-j} + \frac{\partial f}{\partial \varepsilon_{k}}
\bigtriangleup \varepsilon_{k}, \label{eq:StochasticEq21}
\end{equation}
where $\bigtriangleup \mathrm{y}_{k}= \mathrm{y}_k -
\mathrm{y}_{k-1}$, $\bigtriangleup u_{k} = u_k - u_{k-1}$ and
$\bigtriangleup \varepsilon_{k}= \varepsilon_{k} -
\varepsilon_{k-1}$. Note also that without loss of generality we
assumed that $\tau = 0$ in order to simplify notation in the
following discussion. Applying the delay operator, $z^{-1}$ to
both sides of~\eqref{eq:StochasticEq21}, yields
\begin{equation}
\bigg(1 - \sum_{i=1}^n  \frac{\partial f}{\partial
\mathrm{y}_{k-i}} z^{-i} \bigg) \bigtriangleup \mathrm{y}_{k} =
\bigg(\sum_{j=0}^m \frac{\partial f}{\partial u_{k-j}}
z^{-j}\bigg) \bigtriangleup u_{k} + \frac{\partial f}{\partial
\varepsilon_{k}} \bigtriangleup \varepsilon_{k}.
\label{eq:StochasticEq22}
\end{equation}
Given assumption $1$ and that $Q$ is selected such that,
\begin{equation}
{2Q+ \frac{\partial \varphi}{\partial u} \bigg |_{k-1}} \neq 0,
\label{eq:StochasticEq23}
\end{equation}
yields that the control input as calculated
from~\eqref{eq:StochasticEq19} is bounded. This together with
assumption $2$ allows generating the following stochastic model
for the control input
\begin{equation}
u^*_k = u_k + e_{u_k}, \label{eq:StochasticEq24}
\end{equation}
where $u_k$ denotes the control input as estimated from the
controller MDN. In other words, the error as a result of
estimating the control signal from the controller MDN, $e_{u_k}$
is bounded. From~\eqref{eq:StochasticEq19}
and~\eqref{eq:StochasticEq24},
\begin{equation}
\bigtriangleup u_k = B_k \bigtriangleup \chi_k - \bigtriangleup
e_{u_k}, \label{eq:StochasticEq25}
\end{equation}
where we have defined the following row vector
\begin{equation}
B_k = -\frac{\frac{\partial \varphi}{\partial \chi} \bigg
|_{k-1}}{2Q + \frac{\partial \varphi}{\partial u} \bigg |_{k-1}}
\in \mathbb{R}^{1\times(n+m+1)}. \label{eq:StochasticEq26}
\end{equation}
Then, $\bigtriangleup u_k$ can be expressed as
\begin{equation}
\bigtriangleup u_k = \sum_{i=1}^n B_k^i \bigtriangleup
\mathrm{y}_{k-i} + \sum_{j=1}^m B_k^{n+j} \bigtriangleup u_{k-j} +
B_k^{n+m+1} \bigtriangleup \mathrm{y}^d_{k} - \bigtriangleup
e_{u_k}. \label{eq:StochasticEq27}
\end{equation}
Solving~\eqref{eq:StochasticEq27} for $\bigtriangleup u_k$ ,
yields
\begin{equation}
\bigtriangleup u_k = \frac{\sum_{i=1}^n B_k^i z^{-i}
\bigtriangleup \mathrm{y}_{k} + B_k^{n+m+1} \bigtriangleup
\mathrm{y}^d_{k} - \bigtriangleup e_{u_k}}{1-\sum_{j=1}^m
B_k^{n+j} z^{-j}}. \label{eq:StochasticEq28}
\end{equation}
Using~\eqref{eq:StochasticEq28} in~\eqref{eq:StochasticEq22}
yields the following equation for the linearized closed loop
system,
\begin{equation}
N(z^{-1},k) \bigtriangleup \mathrm{y}_{k} = M(z^{-1},k) \bigg
(B_k^{n+m+1} \bigtriangleup \mathrm{y}^d_{k} - \bigtriangleup
e_{u_k} \bigg) + \omega_k \bigtriangleup \varepsilon_{k},
\label{eq:StochasticEq29}
\end{equation}
where
\begin{eqnarray}
N(z^{-1},k) &=& 1 - \sum_{j=1}^m B_k^{n+j} z^{-j} - \sum_{i=1}^n
\frac{\partial f}{\partial \mathrm{y}_{k-i}} z^{-i} + \sum_{i=1}^n
\frac{\partial f}{\partial \mathrm{y}_{k-i}} z^{-i} \sum_{j=1}^m
B_k^{n+j} z^{-j} - \sum_{j=0}^m \frac{\partial f}{\partial
u_{k-j}} z^{-j} \sum_{i=1}^n B_k^i z^{-i}
\label{eq:StochasticEq30} \\
M(z^{-1},k) &=& \sum_{j=0}^m \frac{\partial f}{\partial u_{k-j}}
z^{-j} \label{eq:StochasticEq31} \\
\omega_k &=& \bigg ( 1- \sum_{j=1}^m B_k^{n+j} z^{-j} \bigg)
\frac{\partial f}{\partial \varepsilon_{k}},
\label{eq:StochasticEq32}
\end{eqnarray}
and where $N(z^{-1},k)$ and $M(z^{-1},k)$ are two polynomials that are
related to the system and control structure. Define $\eta_k
=M(z^{-1},k) \bigg (B_k^{n+m+1} \bigtriangleup \mathrm{y}^d_{k} -
\bigtriangleup e_{u_k} \bigg) + \omega_k \bigtriangleup
\varepsilon_{k} $, then the linearized closed loop system
of~\eqref{eq:StochasticEq29} can be rewritten as,
\begin{equation}
N(z^{-1},k) \bigtriangleup \mathrm{y}_{k} = \eta_k.
\label{eq:StochasticEq33}
\end{equation}
From assumption $2$ and the fact that $\mathrm{y}^d_{k}$ and
$\varepsilon_{k}$ are bounded yield that $\eta_k$ is bounded. This
means that the closed loop system is stable if we can guarantee
that $\bigtriangleup \mathrm{y}_{k}$ is bounded. This can be found
out by obtaining the state space representation of $\bigtriangleup
\mathrm{y}_{k}$. For that purpose define,
\begin{eqnarray}
N(z^{-1},k) &=& 1 - \sum_{i=1}^{n+m} a_k^i z^{-i},
\label{eq:StochasticEq34} \\
\mathbb{X}_k &=& [x_k^1, x_k^2, \ldots, x_k^{n+m}]^T
\label{eq:StochasticEq35}
\end{eqnarray}
where $x_k^1 = \bigtriangleup \mathrm{y}_{k-m-n}, x_k^2 =
\bigtriangleup \mathrm{y}_{k-m-n+1}, \ldots, x_k^{n+m} =
\bigtriangleup \mathrm{y}_{k-1}$. Then using the direct programming
method, the state space representation of $\bigtriangleup
\mathrm{y}_{k}$ is given by,
\begin{equation}
\mathbb{X}_{k+1}= \left [ \begin{array}{cccc} 0 &~ 1 &~\ldots &~0
\\ \vdots &~\vdots &~\vdots &~\vdots  \\ 0 &~ 0 &~0 &~1 \\
a^{n+m}_k &~a^{n+m-1}_k &~\ldots &~a^{1}_k
\end{array} \right ] \mathbb{X}_k + \left [ \begin{array}{c} 0 \\ \vdots \\
0 \\ 1 \end{array} \right ] \eta_k, \label{eq:StochasticEq36}
\end{equation}
From~\eqref{eq:StochasticEq36}, it can be seen that the condition
for the local stability of the closed loop system is given by
\begin{equation}
\parallel A_k \parallel < 1,
~~~~~~~~~~~~\text{with}~~~~~~~~~~~A_k = \left [ \begin{array}{cccc} 0 &~ 1 &~\ldots &~0 \\ \vdots
&~\vdots &~\vdots &~\vdots  \\ 0 &~ 0 &~0 &~1 \\ a^{n+m}_k
&~a^{n+m-1}_k &~\ldots &~a^{1}_k
\end{array} \right ]. \label{eq:StochasticEq38}
\end{equation}
\section{Numerical Simulations and Results}
The advantages of the proposed probabilistic control approach are
evaluated on two test stochastic control examples with additive
and multiplicative non--Gaussian random noise. In this section we
provide a comparison of the proposed probabilistic control method
with the conventional indirect adaptive control approach. At the
same time, these results also demonstrate that the inclusion of
models uncertainty significantly enhances the performance of
control systems for the two nonlinear stochastic control examples.

\subsection{Example 1}
\begin{figure*}[htbp]
\centerline{\subfigure[]{\includegraphics[width=1.5in]{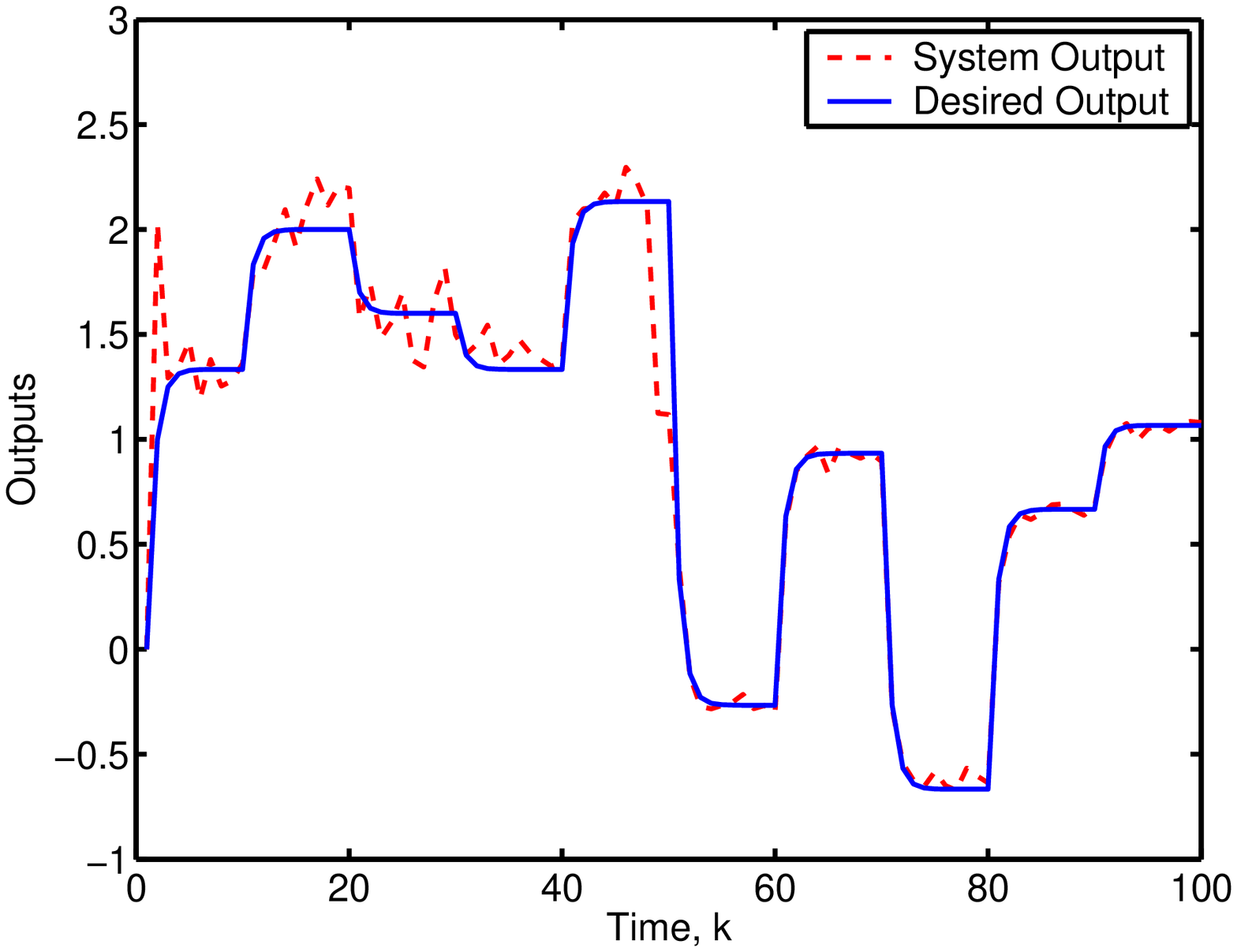}
\label{fig_first_case}} \hfil
\subfigure[]{\includegraphics[width=1.5in]{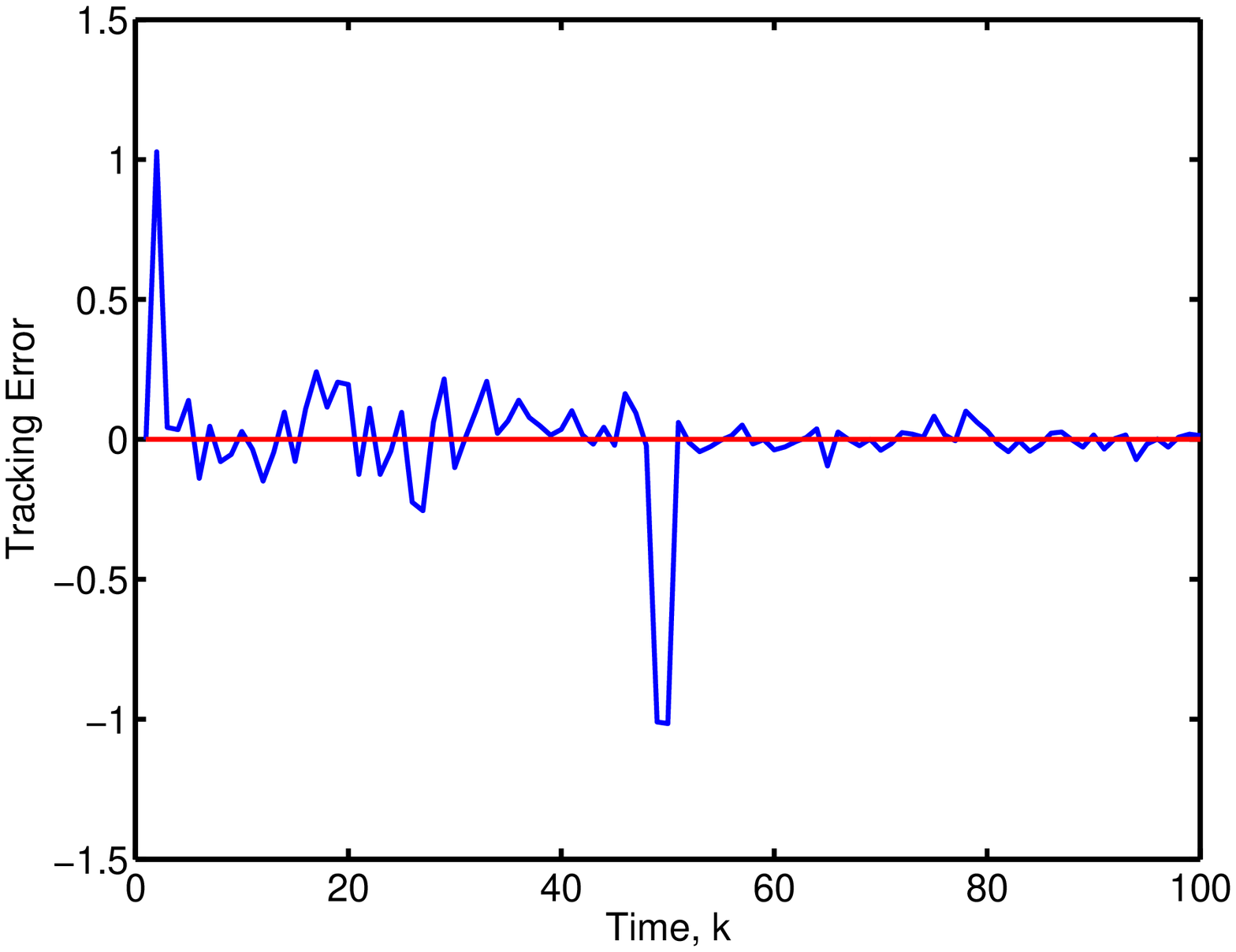}
\label{fig_second_case}} \vspace*{4pt}
{\includegraphics[width=1.5in]{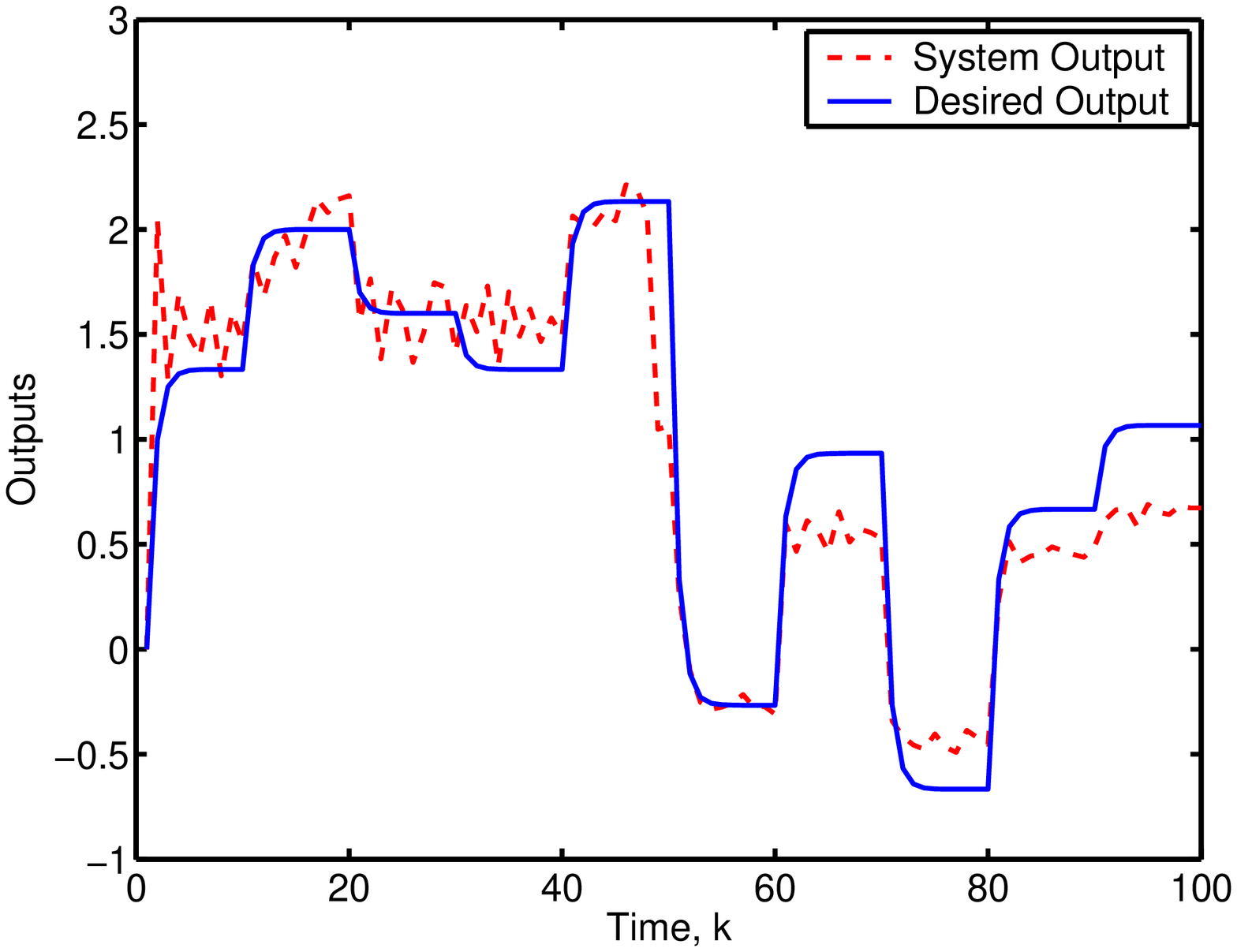}
\label{fig_third_case}} \hfil
\subfigure[]{\includegraphics[width=1.5in]{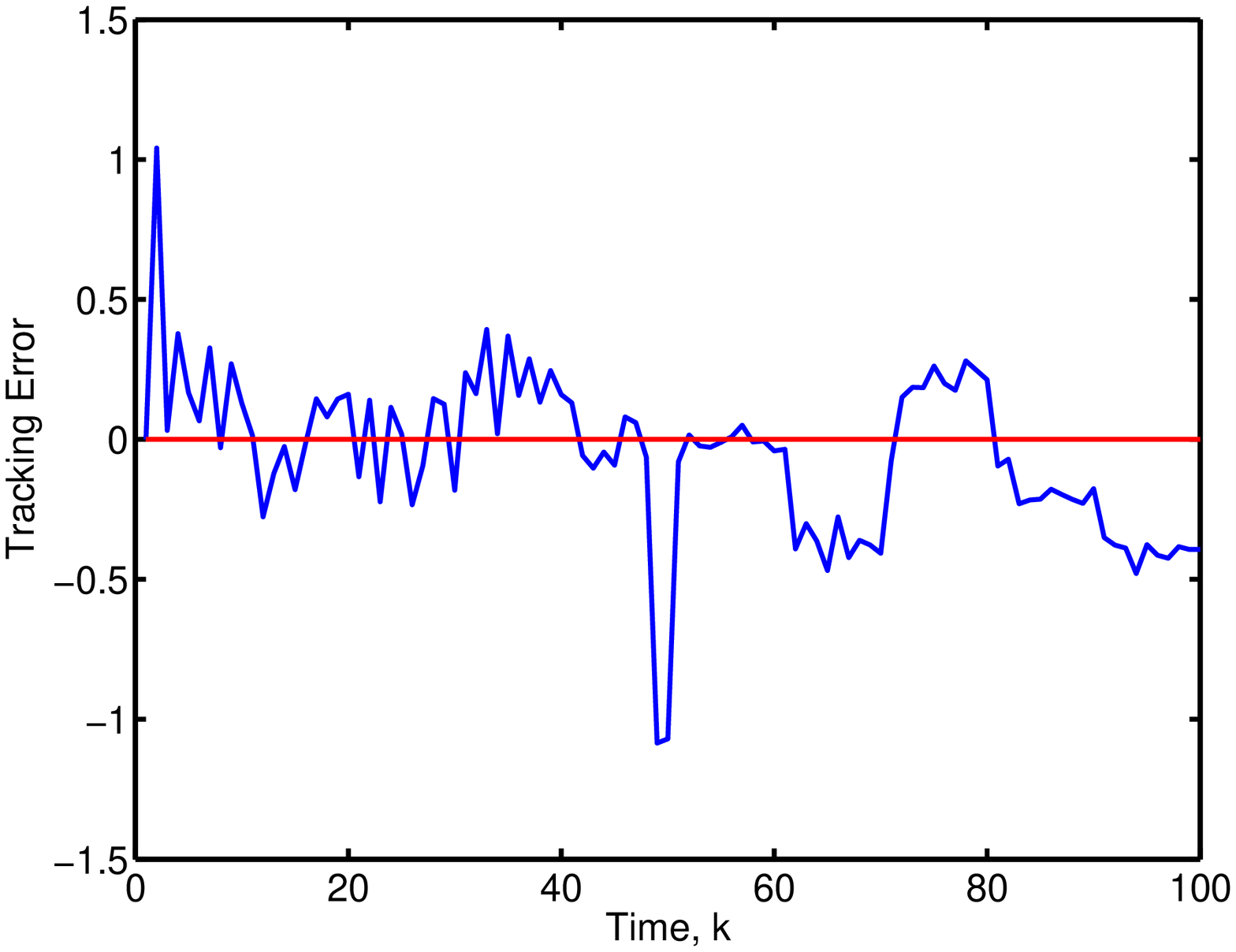}
\label{fig_fourth_case}}} \caption{Control results of nonlinear
stochastic system with additive noise: output and tracking error
(a) the actual and reference model outputs of MDN. (b) tracking error of MDN. (c) the actual and reference model outputs of standard MLPN. (d) tracking error of standard MLPN.} \label{fig:MixFig1}
\end{figure*}

A nonlinear stochastic control problem is
considered here to test the effectiveness of the proposed probabilistic
control approach. The nonlinear stochastic dynamical system is
described by the following equation
\begin{equation}
\mathrm{y}_k = 0.5 - 0.02\mathrm{y}_{k-1}(1.2+\tan^{-1}u_k)+0.2
u_k + \varepsilon_k, \label{eq:StochasticEq39}
\end{equation}
where $\varepsilon_k$ denotes a noise sequence sampled from a
mixture of Gaussians with the following mean, $\mu_{\varepsilon_k}$
and covariance matrix, $\Sigma_{\varepsilon_k}$,
\begin{eqnarray}
\mu_{\varepsilon_k} = \left [ \begin{array}{cc} 1 &~ 0
\end{array} \right ],  &\Sigma_{\varepsilon_k} = \left [
\begin{array}{cc}  0.02 &~ 0\\0 &~0.001  \end{array} \right ].
\nonumber
\end{eqnarray}
The following reference model with input output pairs
$\{r_k,\mathrm{y}^d_k\}$ represents the desired output behavior at
time $k$,
\begin{equation}
\mathrm{y}^d_k = r_k + 0.25 \mathrm{y}^d_{k-1}.
\label{eq:StochasticEq40}
\end{equation}
For comparison purposes two experiments were conducted. In the
first experiment the dynamical behavior of the system was
estimated using a standard multi layer perceptron neural network
(MLPNN) and the controller was estimated using the classical
indirect adaptive control method such that the following error is
minimized,
\begin{equation}
e_k^t=(\mathrm{y}_k - \mathrm{y}_k^d)^2 + Q u_k^2.
\label{eq:StochasticEq41}
\end{equation}
Here the inverse controller was a standard MLPNN as well. In the
second experiment, the conditional distribution of the tracking
error was estimated using the MDN and the
controller was derived from this distribution and estimated using
another MDN as discussed in
Section~\ref{sec:Stat}. In both experiments the scaled
conjugate gradient method is used to update the networks parameters. For fair comparison, both of the MDN and standard neural networks were subjected to the same noise sequence, reference input, and weights $R=0.4, M=1.0$ and
$Q=0.001$ for the variance mean and control input respectively. A rough initialization for the parameters of the standard MLP and
MDNs in both experiments was obtained using an
off line training method. The optimal control law from the
standard MLP control model and the most probable control from the
density network were calculated and forwarded to the system. The
performance of the two controllers is shown in
Figure~\ref{fig:MixFig1}. This figure shows the system output
as a result of the MDN control superimposed on the desired
model output over the whole control range. Clearly, the system
output is able to track adequately the desired output. However, the
standard NN output is struggling to track the desired
output. This shows that although standard NNs are
normally suitable for deriving the optimal control law and
achieving a good tracking performance, they are inadequate for
stochastic systems affected by general random inputs. Looking at
the tracking error of the two methods it can be seen that on
average the tracking error of the MDN controller is zero. However,
the average tracking error of the standard NN is drifted away from zero.
\subsection{Example 2}
\begin{figure*}[htbp]
\centerline{\subfigure[]{\includegraphics[width=1.5in]{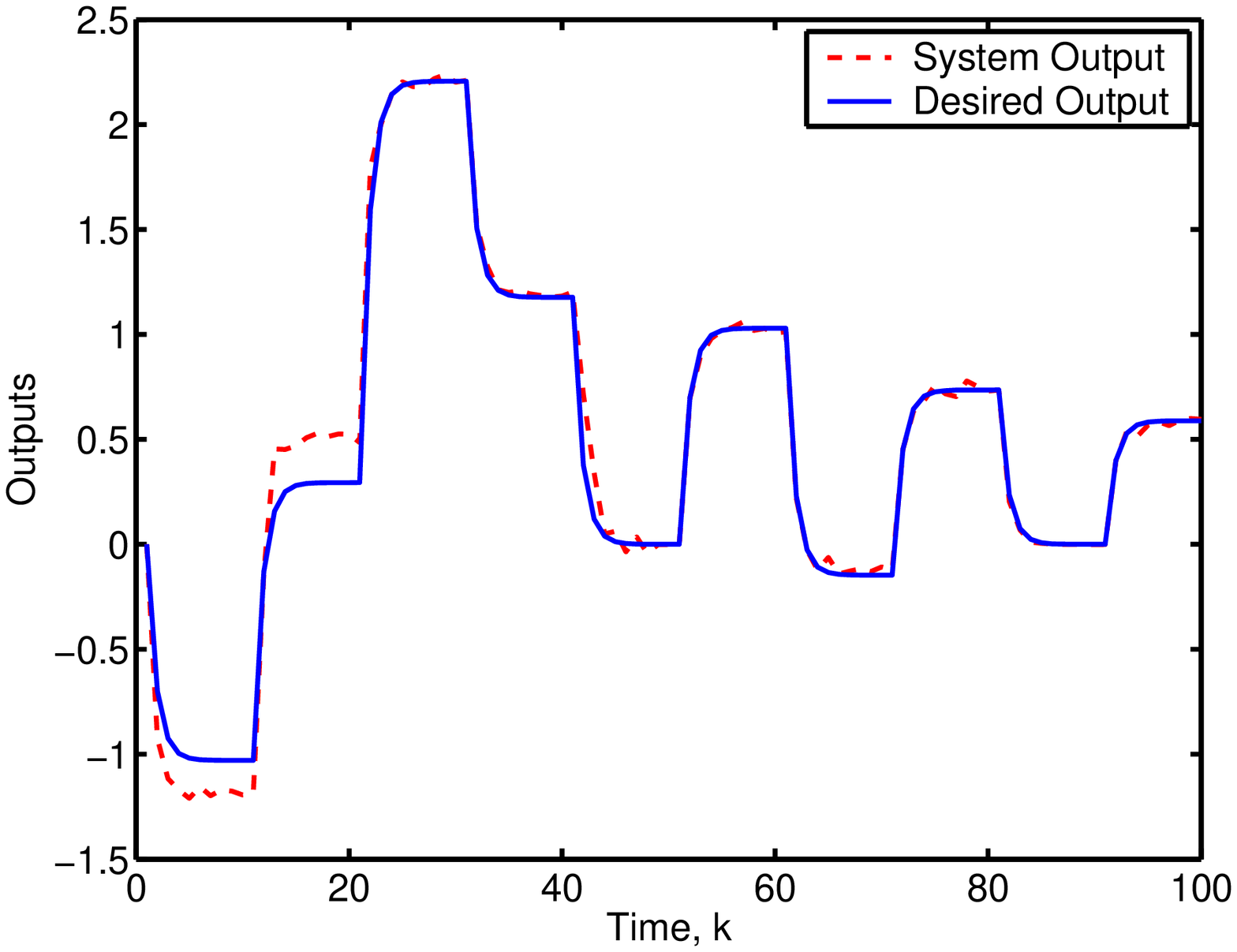}
\label{fig_first_case}} \hfil
\subfigure[]{\includegraphics[width=1.5in]{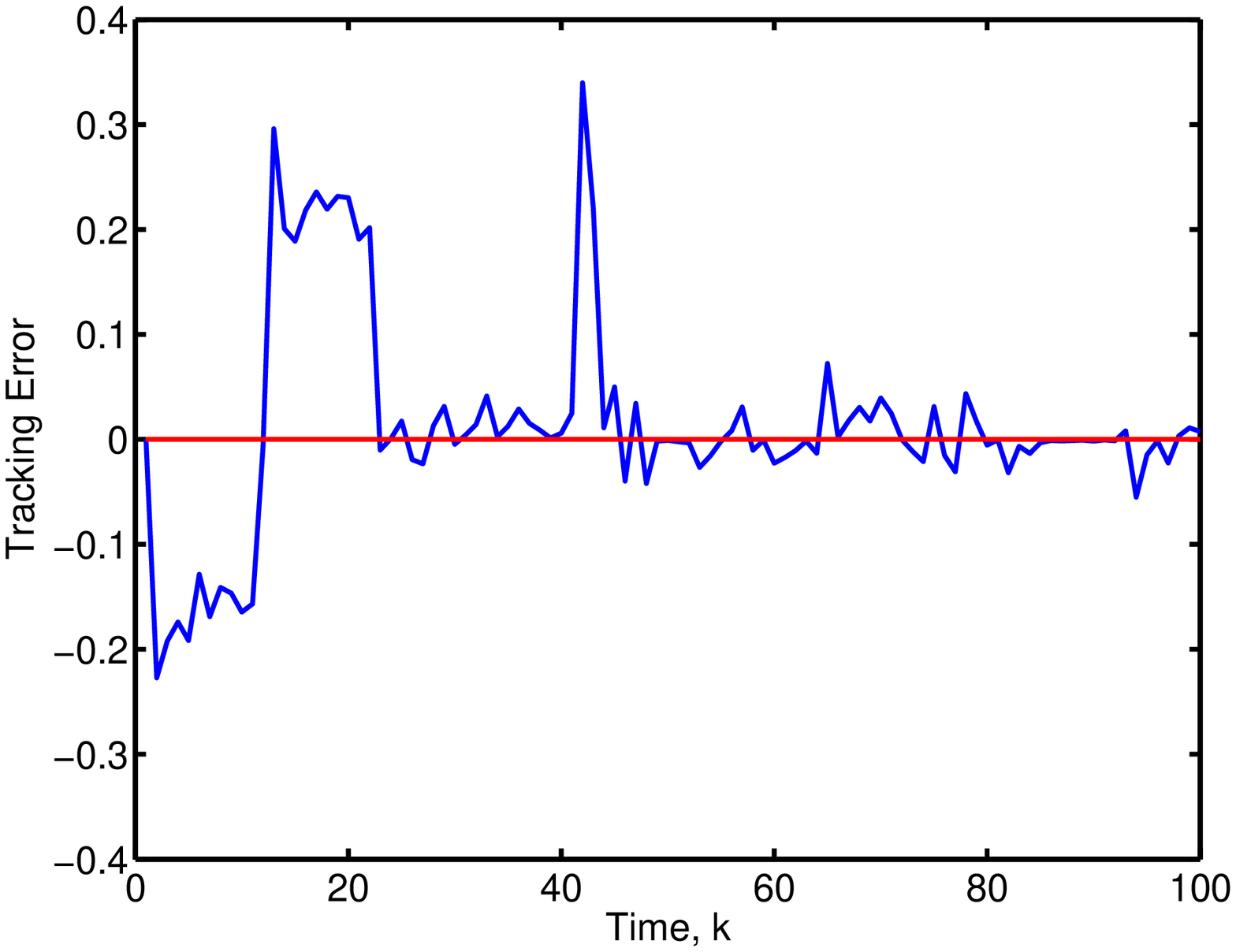}
\label{fig_second_case}} \hfil {\includegraphics[width=1.5in]{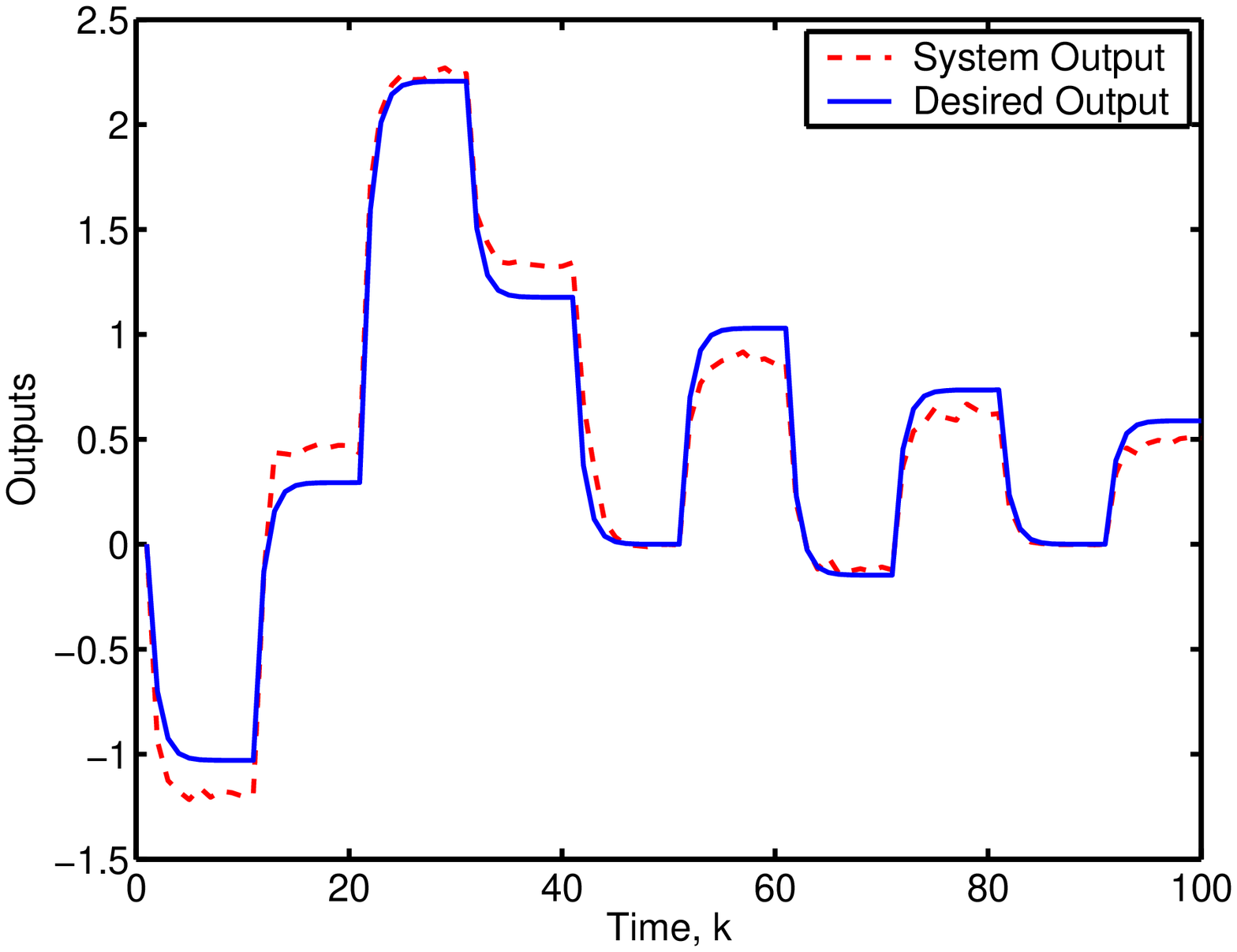}
\label{fig_third_case}} \hfil
\subfigure[]{\includegraphics[width=1.5in]{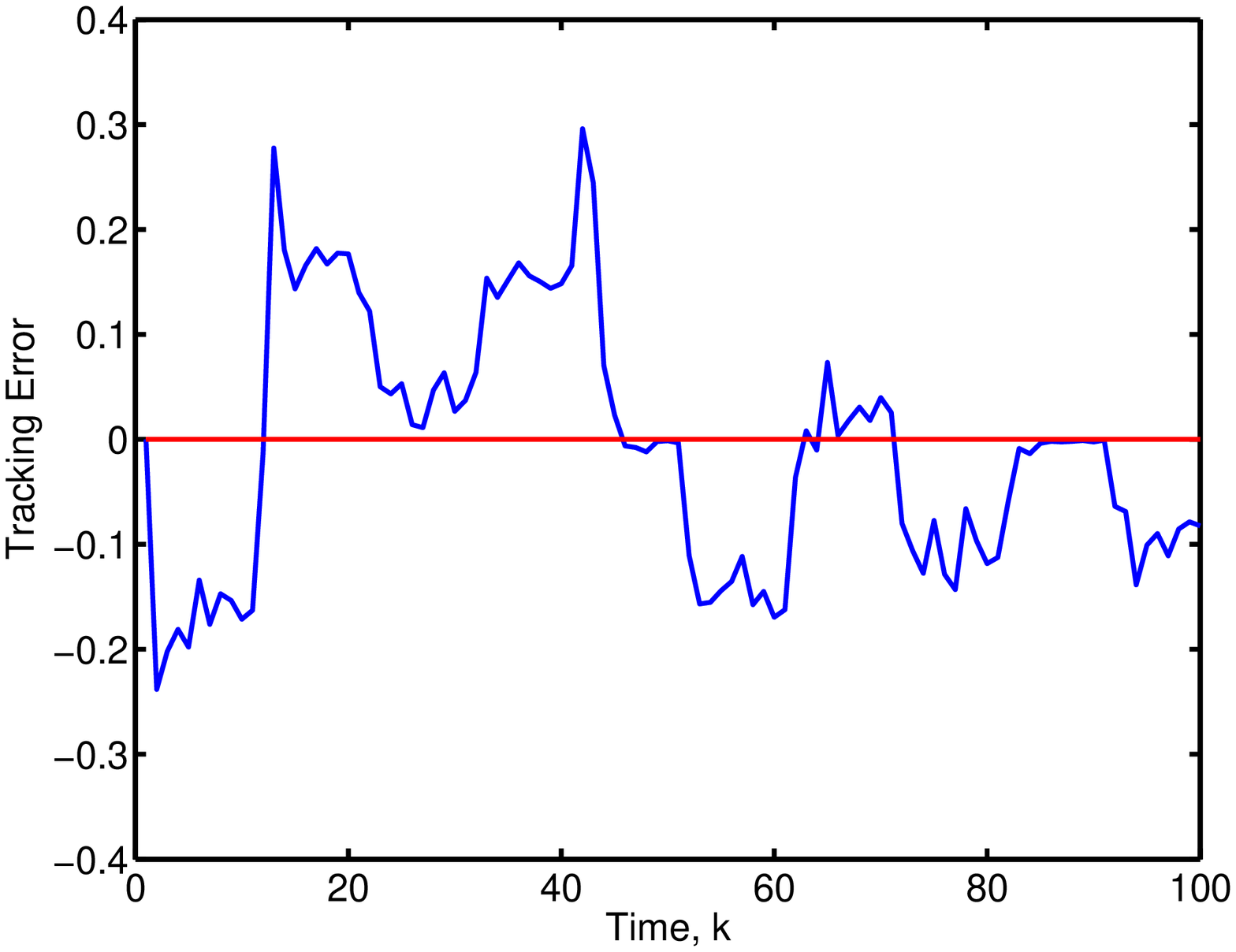}
\label{fig_fourth_case}}} \caption{Control results of nonlinear
stochastic system with multiplicative noise: output and tracking
error (a) the actual and reference model outputs of MDN. (b) tracking error of MDN. (c) the actual and reference model outputs of standard MLPN. (d) tracking error of standard MLPN.} \label{fig:MixFig2}
\end{figure*}

In this section we present a comparison of the proposed
probabilistic control approach with the indirect adaptive control
approach on a stochastic nonlinear control problem with
multiplicative noise. The stochastic nonlinear dynamical system is
described by the following difference equation:
\begin{equation}
\mathrm{y}_k = \frac{29}{40} \varepsilon_k \sin \bigg (\frac{16
u_{k-1} + 8 \mathrm{y}_{k-1}}{3 + 4 u_{k-1}^2 + 4
\mathrm{y}_{k-1}^2} \bigg) + \frac{2}{10} (u_{k-1} +
\mathrm{y}_{k-1}), \label{eq:StochasticEq42}
\end{equation}
where $\varepsilon_k$ denotes a noise sequence sampled from a
mixture of Gaussians with the following mean, $\mu_{\varepsilon_k}$
and covariance matrix, $\Sigma_{\varepsilon_k}$:
\begin{eqnarray}
\mu_{\varepsilon_k} = \left [ \begin{array}{cc} 0.5 &~ 0
\end{array} \right ],  &\Sigma_{\varepsilon_k} = \left [
\begin{array}{cc} 0.002 &~ 0 \\0 &~ 0.001 \end{array} \right ].
\nonumber
\end{eqnarray}
The particular system considered here was used in~\cite{Singla07},
but with an additive Gaussian white noise sequence rather than
multiplicative non--Gaussian noise. The following reference model
with input output pairs $\{r_k,\mathrm{y}^d_k\}$ represents the
desired output behavior at time $k$
\begin{equation}
\mathrm{y}^d_k = r_{k-1} + 0.32 \mathrm{y}^d_{k-1}.
\label{eq:StochasticEq43}
\end{equation}
Two experiments were conducted. In the first experiment
two standard multi layer perceptron NNs were used to
estimate the forward dynamics and the inverse controller of the
system. In the second experiment, two MDNs were used to provide estimates for the conditional distributions of the tracking error and the inverse controller. Here also the scaled conjugate gradient method is used to update the networks parameters. Both of the MDN and standard NNs were subjected to the same noise
sequence, reference input, and weights $R=0.25, M=1.0$ and
$Q=0.01$ for the variance mean and control input respectively. A rough initialization for the parameters of the standard MLP and
MDNs in both experiments was obtained using
off line training methods. The optimal control law from the
standard MLP control model and the most probable control from the
density network were calculated and forwarded to the system. The
performance of the two controllers is shown in
Figure~\ref{fig:MixFig2}. The result of this example is
consistence with that obtained in the first example. The system
output obtained from the MDN controller is superimposed on the
desired model output over all the control range and the tracking
error on average is equal to zero. However the standard neural
network output is struggling to track the desired output and its
average tracking error is drifting away from zero.
\section{Conclusion}
In this paper a new framework has been applied to the design of
controllers which encompasses uncertainty, multimodality,
hysteresis, and arbitrary density functions of forward models and
inverse controllers. It is for the general class of stochastic
nonlinear control problems, where the dynamics are nonlinear
functions of the control and the state. The proposed framework
considers functional uncertainty by estimating the probabilistic
models of the system dynamics and the controller that are
dependent on the state and the control input. Using mixture
density networks from the neural network field, a control input is
formulated which minimizes uncertainty of the closed loop system.
Furthermore, the local stability condition of the closed loop
system has been established. Global stability and closed-loop
performance are topics of future research.

In contrast to traditional control methods, the derived controller in this paper is not constrained by the probability density
function of the random input that affects the dynamics of the
system. No assumptions are made about the invertibility of the
system output with respect to the random input or even the
invertibility of the control input with respect to the system
output. Simulation examples with additive and multiplicative
random inputs are used to illustrate the proposed controller and
encouraging results have been obtained.
\\
\\
\textbf{Acknowledgment}: This work has been carried out during sabbatical leave granted to the author Randa Herzallah from Al-Balqa' Applied University (BAU) during the academic year $2010/2011$


\end{document}